\documentclass[11pt,a4paper]{article}

\usepackage[T1]{fontenc}
\usepackage[utf8]{inputenc}
\usepackage{lmodern}
\usepackage[margin=25mm]{geometry}
\usepackage{amsmath,amssymb,amsthm}
\usepackage{bm}
\usepackage{graphicx}
\usepackage{subcaption}
\usepackage{booktabs}
\usepackage{float}
\usepackage[section]{placeins}
\usepackage{authblk}
\usepackage[numbers,sort&compress]{natbib}
\usepackage{microtype}
\usepackage{xcolor}
\usepackage[hidelinks]{hyperref}

\hypersetup{
  pdftitle={Reconstructing hidden attractors with reservoir computing},
  pdfauthor={Mohammad Zaman, Emel Khan, Thomas Geert de Jong}
}

\newtheorem{theorem}{Theorem}
\newtheorem{lemma}[theorem]{Lemma}

\title{Reconstructing hidden attractors with reservoir computing}
\author[1]{Mohammad Zaman}
\author[2]{Emel Khan}
\author[3,4]{Thomas Geert de Jong\thanks{Email: \href{mailto:t.g.de.jong.math@gmail.com}{\nolinkurl{t.g.de.jong.math@gmail.com}}}}
\affil[1]{Government Post Graduate College, Charsadda 24420, Khyber Pakhtunkhwa, Pakistan}
\affil[2]{Institute of Numerical Sciences, Kohat University of Science and Technology, Kohat 26000, Khyber Pakhtunkhwa, Pakistan}
\affil[3]{Temple University, Japan Campus, 1-14-29 Taishido, Setagaya-ku, Tokyo 154-0004, Japan}
\affil[4]{School of Mathematical Sciences, Hebei Normal University, 20 South Erhuan East Road, Shijiazhuang 050024, Hebei, China}
\date{}

\begin{document}
\maketitle

\begin{abstract}
Hidden attractors are difficult to locate as their basins of attraction do not intersect with any small neighborhood of an equilibrium point. We examine the possibility of reconstructing the hidden
chaotic attractor from time series via reservoir
computing, where the reservoir is built by a network of driven Kuramoto oscillators.  We compare all-to-all coupling, nearest-neighbour ring coupling, and a ring with one modified oscillator equation which ensures that any attractor will be hidden. The numerical results suggest that ring coupling supports accurate reconstruction over a wider range of coupling and forcing strengths than all-to-all coupling. Searches near reservoir equilibria nevertheless show why agreement in short-term prediction and long-term statistics does not establish preservation of hiddenness. Numerics for the modified ring show that it has the highest performance which coincides with our intuition that reconstructing the global dynamics would yield a more faithful attractor reconstruction. We also show that this construction precludes global forgetting of initial conditions for constant inputs, although local conditional contraction along a nonconstant signal remains possible. The distinction between reconstructing attractor dynamics and preserving basin structure provides a basis for evaluating and designing reservoir models of hidden attractors.
\end{abstract}


\section{Introduction}\label{sec:introduction}

An attractor is hidden when its basin cannot be reached from arbitrarily small neighbourhoods of equilibria. This separates hidden attractors from self-excited attractors, which can be located by following trajectories initiated near an unstable equilibrium \citep{leonov2012hidden,leonov2013hidden,dudkowski2016hidden}. Systems without equilibria provide a particularly direct setting: any attracting oscillation is hidden under this classification. The distinction concerns the surrounding phase space, rather than only the motion observed on the attractor.

Reservoir computing learns dynamical behaviour by driving a fixed nonlinear system with data and fitting a linear readout \citep{jaeger2001echo,jaeger2004harnessing,lukovsevivcius2009reservoir}. Successful models can support short-term forecasting and reproduce long-term chaotic statistics \citep{pathak2017using,pathak2018model,lu2018attractor}. Physical oscillator networks provide one implementation of this framework \citep{nakajima2020physical,de2025harnessing}. Here we use driven Kuramoto oscillators, whose closed-loop equilibria can be studied explicitly.

Preserving the hidden character of an attractor is a stronger requirement than reproducing its trajectory or invariant measure. Even exact correspondence between dynamics restricted to two attractors would not specify the equilibria and basins outside them. Related studies have examined inference of unseen attractors \citep{rohm2021model,gauthier2022learning}, extrapolation from local dynamical examples \citep{kim2021teaching}, and the difficulty of predicting basins with reservoir models \citep{zhang2023catch}. These results motivate a separate question: when a reservoir reconstructs a hidden attractor, is the reconstructed attractor hidden in the reservoir's own state space?

We study this question for the NE6 system \citep{jafari2013elementary}. The numerical comparison concerns three reservoir architectures: all-to-all coupling, a nearest-neighbour ring, and a ring in which one phase equation is changed so that the resulting reservoir has no equilibria. We distinguish prediction error, long-term output statistics, local conditional contraction, and equilibrium--basin relationships. We also establish two analytical results: the modified ring has no equilibria for any trained readout, and this absence obstructs global initial-state forgetting for constant input.

The paper is structured as follows. Section \ref{Hidden attractors} defines the hidden attractor and Section \ref{The NE6 Chaotic system} the NE6 system. Section \ref{general RC} provides the overview of the general scheme of reservoir computing (RC), and Section \ref{Kuramoto oscillators} the Kuramoto reservoir with three coupling topologies. Section \ref{metrics} provides the definitions of the metrics used to characterize the performance of the systems and the classification of the reconstructed attractors into hidden and self-excited attractors. Section~\ref{Results} presents the results for the all-to-all, ring and
modified ring reservoirs, and Section~\ref{discussion} discusses them,
the limitations and future work.

\section{Hidden attractors}\label{Hidden attractors}

Hidden attractors are attractors whose basin of attraction does not intersect with the neighbourhoods of any equilibrium point. The idea of hidden attractors was first coined by Leonov and Kuznetsov \cite{leonov2013hidden,leonov2011localization,leonov2012hidden} when they found a hidden attractor in the chaotic Chua's circuit. Unlike self-excited attractors that have basins of attraction connected with the unstable equilibria and therefore can be easily found using trajectories starting from the neighborhood of equilibria, hidden attractors cannot be found numerically using localization methods based on local linearization.  The localization and visualization
of hidden attractors therefore constitute a challenging computational problem, requiring carefully selected initial conditions and specialized search algorithms. Hidden attractors have been reported in a wide range of nonlinear systems, including mechanical oscillators, electronic circuits,  memristive systems, phase-locked loops, and power systems \cite{kuznetsov2020theory}.

The existence of hidden attractors is especially important in the context of multistable systems, wherein there might be several coexisting attractors having basins of attraction that are hard to discern. Hence, there is always the possibility that the system will suddenly switch from one desired attractor to another unknown or undesired attractor, bringing about a sudden change in the behavior of the system and possibly causing some adverse effect \cite{scheffer2009early,dudkowski2016hidden}. Therefore, it is of prime importance that all the coexisting attractors, including those which are hidden, be identified.

\section{The NE6 hidden attractor}\label{The NE6 Chaotic system}

In this work, we study the NE6 chaotic system proposed by
Jafari \emph{et al.}~\cite{jafari2013simple} as a typical example of a three-dimensional
chaotic system without any equilibria. The NE6 system belongs to a
class of elementary quadratic chaotic flows and reveals hidden chaotic
behavior, even with its simple mathematical representation. Contrary to the
classical chaotic systems like the Lorenz~\cite{lorenz2017deterministic} and R\"ossler~\cite{rossler1976equation} systems, which have equilibria,
the NE6 system does not have any equilibria.

The governing equations of the NE6 system are given by

\begin{equation}
\begin{aligned}
\dot{x} &= y,\\
\dot{y} &= z,\\
\dot{z} &= -y-xz-yz-a,
\end{aligned}
\label{eq:NE6}
\end{equation}

where $a>0$ is the control parameter. Throughout this study, we set
$a=0.75$ and use the initial condition $(x_0,y_0,z_0)=(0,3,-0.1)$,
which generates chaotic dynamics. For these parameter values, the
Lyapunov spectrum is $(0.0280,\,0,\,-3.4341)$, confirming the chaotic
nature of the attractor. The dynamical behavior of the NE6 system can be characterized by varying the parameter $a$. As $a$ is varied, the system undergoes a sequence of period-doubling bifurcations, leading from periodic doubling bifurcation to chaotic dynamics \cite{jafari2013elementary}.

\begin{figure}[tbp]
\centering

\begin{subfigure}[b]{0.49\columnwidth}
    \centering
    \includegraphics[width=\linewidth]{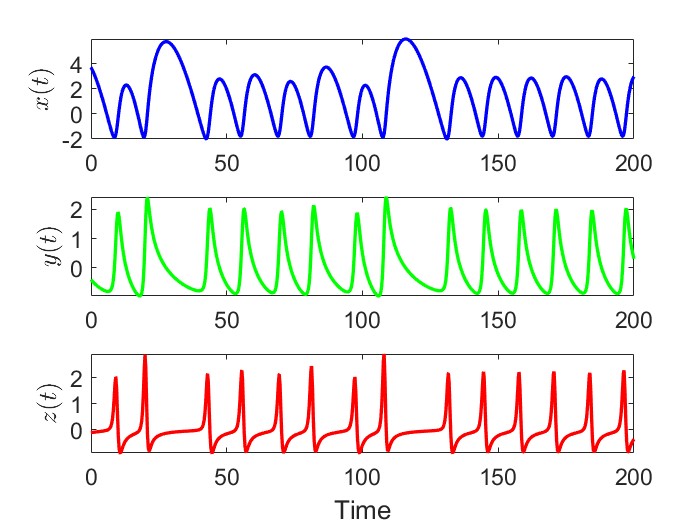}
\end{subfigure}
\hfill
\begin{subfigure}[b]{0.49\columnwidth}
    \centering
    \includegraphics[width=\linewidth]{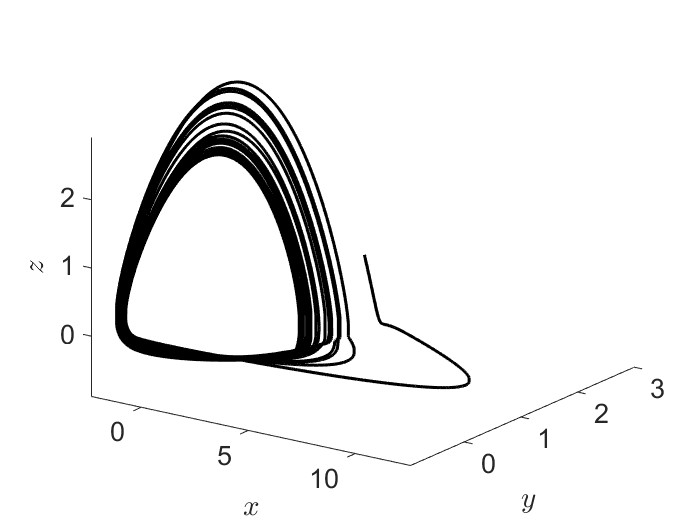}
\end{subfigure}
\caption{(a) Time series of the state variables $x(t)$, $y(t)$, and $z(t)$ of the NE6 chaotic system. (b) Phase-space trajectory of the corresponding hidden chaotic attractor for $a=0.75$ and the initial condition $(x_0,y_0,z_0)=(0,3,-0.1)$.}
\label{fig:NE6}
\end{figure}

The time series and the phase space trajectory of the NE6 attractor are shown in Figure~\ref{fig:NE6}. The output multivariate time series $u(t)=[x(t),y(t),z(t)]^{T}$ is then utilized for training and autonomous predictions in the suggested Kuramoto RC.

\section{Reservoir Computing: General Setup}
\label{general RC}

In this section, we briefly describe the continuous framework of RC~~\cite{lu2018attractor}.

Let $u:\mathbb{R}\to\mathcal{U}$
be an input signal, where $\mathcal{U}$ is an $m$-dimensional manifold
($\mathcal{U}=\mathbb{R}^3$ for NE6).

The reservoir is evolved using the governing equations given by
\begin{equation}
  \frac{dr}{dt}=R(r,u),
  \label{eq:reservoir}
\end{equation}
where the state $r(t)$ resides on an $N$-dimensional manifold $\mathcal{R}$,
the reservoir state space.

Building and testing a reservoir requires three steps:

\begin{enumerate}

\item \emph{Wipe-out.}
Integrate \eqref{eq:reservoir} with the prescribed input $u$ from an
initial condition $r(-T_{\mathrm{wipe}})=r_0$ over
$[-T_{\mathrm{wipe}},0]$. This initial integration reduces the
influence of the initial reservoir state before training begins.
The resulting state $r(0)$ initializes the training phase.

\item \emph{Training.}
Continue integrating the driven reservoir over
$[0,T_{\mathrm{train}}]$.
For an input $u(t)\in\mathbb{R}^{m}$, choose a fixed feature map
$h:\mathcal{R}\to\mathbb{R}^{d_h}$ and approximate the input by
\begin{equation}
  u(t)\approx W_{\mathrm{out}}h(r(t)),
  \qquad
  W_{\mathrm{out}}\in\mathbb{R}^{m\times d_h}.
\end{equation}
Only $W_{\mathrm{out}}$ is trained; the reservoir dynamics and
feature map remain fixed.

Choose sampling times
$0\le t_1<\cdots<t_M<T_{\mathrm{train}}$ and construct the feature
matrix and target matrix
\begin{equation}
\begin{aligned}
  H &=
  \begin{bmatrix}
    h(r(t_1)) & \cdots & h(r(t_M))
  \end{bmatrix}
  \in\mathbb{R}^{d_h\times M},\\
  Y &=
  \begin{bmatrix}
    u(t_1) & \cdots & u(t_M)
  \end{bmatrix}
  \in\mathbb{R}^{m\times M}.
\end{aligned}
\label{eq:training_matrices}
\end{equation}
Each feature vector is paired with the input at the same time,
so that the readout learns a reconstruction of
the driving signal.

The readout weights are obtained by ridge regression:
\begin{equation}
  W_{\mathrm{out}}
  =
  \underset{W\in\mathbb{R}^{m\times d_h}}
  {\operatorname{argmin}}
  \left\{
    \|WH-Y\|_{\mathrm{F}}^{2}
    +\alpha\|W\|_{\mathrm{F}}^{2}
  \right\},
  \qquad \alpha>0,
  \label{eq:readout_regression}
\end{equation}
where $\|\cdot\|_{\mathrm{F}}$ denotes the Frobenius norm.
The first term sums the squared reconstruction errors over all
training samples and input components. The second penalizes
large readout coefficients, improving the conditioning of the
regression when reservoir features are strongly correlated.

We have that
\begin{equation}
  W_{\mathrm{out}}
  =
  YH^{T}(HH^{T}+\alpha I_{d_h})^{-1}.
  \label{eq:readout_solution}
\end{equation}

\item \emph{Testing.}
Fix the trained readout weights and replace the external input
in \eqref{eq:reservoir} by the reconstructed output:
\begin{equation}
\begin{aligned}
  \frac{d\hat r}{dt}
  &=R\bigl(\hat r,W_{\mathrm{out}}h(\hat r)\bigr),\\
  \hat r(T_{\mathrm{train}})
  &=r(T_{\mathrm{train}}).
\end{aligned}
\label{eq:autonomous}
\end{equation}
Integrate this autonomous system over
$[T_{\mathrm{train}},T_{\mathrm{train}}+T_{\mathrm{test}})$.
Its reconstructed signal is
\begin{equation}
  \hat u(t)=W_{\mathrm{out}}h(\hat r(t)).
\end{equation}
The reference signal $u(t)$ is used only to evaluate prediction
accuracy during this phase.

\end{enumerate}

\begin{figure}[tbp]
    \centering \includegraphics[width=0.8\columnwidth]{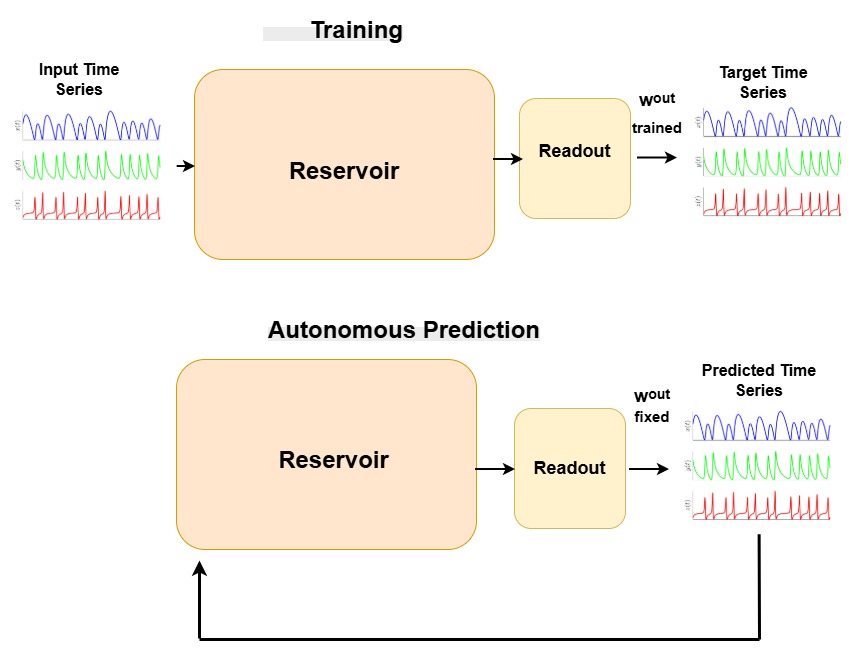}
    \caption{Schematic illustration of RC framework.}
    \label{fig:RC schematic}
\end{figure}

The readout has been trained to replicate the input,
substituting $u$ with $W_{\mathrm{out}}h(r)$ completes the loop. While
the input signal $u(t)$ is used during training, in autonomous
prediction mode it is removed, and the output of the
read-out is recursively fed back into the reservoir
(Fig.~\ref{fig:RC schematic}). In this paper, the vector field $R$ is implemented
by a network of Kuramoto oscillators, as we will see in the next section.

\section{Reservoir computing by coupled oscillators}
\label{Kuramoto oscillators}

The omnipresent computing framework introduced in~\cite{de2025harnessing,de2025designing} states that generic coupled oscillators can be used as a computational resource using the reservoir computing framework. We will use a network of $N$ coupled driven Kuramoto oscillators \cite{kuramoto1984chemical, acebron2005kuramoto}.

\subsection{Driven Kuramoto RC}

The state of the reservoir is defined by the phase of the oscillators
\begin{equation}
r(t)
=
\bm{\theta}(t)
=
\begin{bmatrix}
\theta_1(t) &
\theta_2(t) &
\cdots &
\theta_N(t)
\end{bmatrix}^{T},
\label{eq:reservoir_state}
\end{equation}
where $\theta_i(t)$ denotes the phase of the $i$th oscillator,
$i=1,\ldots,N$. The reservoir manifold is the
$N$-torus
$\mathbb{T}^N:=\underbrace{S\times S\times\cdots\times S}_{N\ \mathrm{times}},
\quad S:=\mathbb{R}/2\pi\mathbb{Z}$, where
i.e., $\bm{\theta}(t)\in\mathbb{T}^N$.

Let $v(t)=(x(t),y(t),z(t))^T$ denote the NE6 state. The
reservoir is driven by the scaled signal
\begin{equation}
  u(t)=\frac{v(t)}{4}+\frac12\bm{1},
  \label{eq:input_scaling}
\end{equation}
where $\bm{1}\in\mathbb{R}^{3}$ is the vector of ones.
Thus, the input dimension is $m=3$.

For the ordinary Kuramoto reservoir, the oscillator phases
satisfy
\begin{equation}
\begin{aligned}
  \frac{d\theta_i}{dt}
  &=
  \omega_i
  +\frac{K}{d}\sum_{j=1}^{N}A_{ij}
       \sin(\theta_j-\theta_i)\\
  &\quad
  +F\sin\bigl(cu_{q_i}(t)-\theta_i\bigr),
  \qquad i=1,\ldots,N.
\end{aligned}
\label{eq:kuramoto}
\end{equation}
Here $\omega_i$ is the natural angular frequency of oscillator
$i$, $K\ge0$ is the coupling strength, and $F\ge0$ is the
forcing amplitude. The adjacency matrix
$A_{ij}\in\{0,1\}$ specifies the coupling graph, and $d$
sets its normalization. Each oscillator receives one input
component, selected by a fixed index
$q_i\in\{1,\ldots,m\}$, while $c$ controls the input gain. The input-assignment vector is
\begin{equation}
  \bm{q}=(q_1,\ldots,q_N)^T
  \in\{1,\ldots,m\}^{N},
\end{equation}
where $q_i=k$ means that oscillator $i$ is driven by the
input component $u_k(t)$. The frequencies, input assignments, input gain, assignment vector are held
fixed throughout the parameter sweep.

We compare three reservoir architectures:
\begin{itemize}

\item \emph{All-to-all coupling.}
Here $A_{ij}=1$ for $i\ne j$, $A_{ii}=0$, and $d=N$.
The coupling term is therefore
\begin{equation}
  \frac{K}{N}\sum_{j=1}^{N}
  \sin(\theta_j-\theta_i),
\end{equation}
where including $j=i$ does not change the sum.

\item \emph{Ring coupling.}
Each oscillator is connected to its two nearest neighbours,
with periodic indexing:
$A_{ij}=1$ when $j=i\pm1\pmod N$ and $A_{ij}=0$
otherwise. We take $d=2$, giving
\begin{equation}
\begin{aligned}
  \frac{K}{2}\sum_j A_{ij}\sin(\theta_j-\theta_i)
  &=
  \frac{K}{2}\sin(\theta_{i-1}-\theta_i)\\
  &\quad+
  \frac{K}{2}\sin(\theta_{i+1}-\theta_i).
\end{aligned}
\label{eq:ring_coupling}
\end{equation}
We assume $N\ge3$. The total incident coupling weight is
$K$ for the ring and $K(N-1)/N$ for the all-to-all network.

\item \emph{Modified ring.}
The adjacency matrix and coupling normalization are unchanged,
but the equation for the last oscillator is replaced.
Let $b_i(\bm{\theta},u)$ denote the right-hand side of
\eqref{eq:kuramoto} with ring coupling. The modified dynamics
are
\begin{equation}
\begin{aligned}
  \dot{\theta}_i
  &=b_i(\bm{\theta},u),
  \qquad i=1,\ldots,N-1,\\
  \dot{\theta}_N
  &=\omega_N+\frac{4\beta}{F}
     \bigl(b_N(\bm{\theta},u)-\omega_N\bigr)
     s(\bm{\theta},u),
\end{aligned}
\label{eq:modified_ring}
\end{equation}
where
\begin{equation}
  s(\bm{\theta},u)
  =
  \frac{
    b_{N-1}(\bm{\theta},u)+b_1(\bm{\theta},u)
  }{2}.
\end{equation}
The implementation uses $\beta=0.01$ and requires $F>0$
and $\omega_N\ne0$.
The modification uses the ordinary velocities of the last
oscillator's two neighbours. Consequently, it changes the
node dynamics rather than the coupling topology.
At a hypothetical equilibrium, the first $N-1$ equations
would give $b_1=b_{N-1}=0$, leaving
$\dot{\theta}_N=\omega_N\ne0$. Thus the modified ring has
no equilibria.

\end{itemize}

\subsection{Numerical setup}

The NE6 time series is generated and sampled
at a fixed interval before it is used as the input $u(t)$. The reservoir
equations are integrated numerically with a fourth-order Runge-Kutta scheme,
and the read-out $W_{\mathrm{out}}$ is computed by linear regression on the
training interval. The numerical parameters, such as the
number of oscillators and the durations of the wipe-out, training and testing
phases, are listed in Table~\ref{tab:parameters}. The parameter plane $(K,F)$ is
sampled on a regular grid of $420$ values, with the ranges shown on the axes
of Figs.~\ref{fig:naive_all_to_all}, \ref{fig:RC_Ring} and \ref{fig:RC_ModRing}, and the
three reservoirs are compared on the same grid and the same NE6 data. Due to
the structure of the vector fields, hardware acceleration was applied to the
simulations using the techniques of \cite{de2023mlncp}.

\begin{table}[tbp]
\centering
\small
\caption{Notation used in the Kuramoto reservoir models.}
\label{tab:parameters}
\begin{tabular}{
  @{}p{0.27\linewidth}
     p{0.67\linewidth}@{}
}
\toprule
Symbol & Description \\
\midrule
$v(t)$
& NE6 state \\

$u(t)$
& Scaled driving signal and training target \\

$m$
& Dimension of the input and reconstructed output \\

$a$
& NE6 control parameter \\

$N$
& Number of reservoir oscillators \\

$\theta_i(t)$
& Phase of oscillator $i$ \\

$\omega_i$
& Natural angular frequency of oscillator $i$ \\

$K$
& Coupling strength \\

$A_{ij}$
& Adjacency matrix of the coupling graph \\

$d$
& Coupling normalization factor \\

$F$
& Forcing amplitude \\

$c$
& Input gain \\

$\bm{q}$
& Fixed vector assigning input components to oscillators \\

$h(\bm{\theta})$
& Periodic feature map containing a constant, sines and squared sines \\

$d_h$
& Number of readout features \\

$W_{\mathrm{out}}$
& Trained linear readout matrix, including the intercept \\

$\alpha$
& Ridge-regression regularization parameter \\

$\beta$
& Gain in the modified-ring equation \\

$T_{\mathrm{wipe}}$
& Duration of the initial wipe-out phase \\

$T_{\mathrm{train}}$
& Duration of the driven training phase \\

$T_{\mathrm{test}}$
& Duration of the autonomous testing phase \\
\bottomrule
\end{tabular}
\end{table}

\section{Performance metrics}
\label{metrics}
We assess prediction accuracy, chaotic dynamics, conditional
stability and distributional agreement using four metrics.

\subsection{Leading Lyapunov exponent}
We estimate the leading Lyapunov exponent (LLE) from the
standardized multivariate prediction $\hat u(t)$ using a
nearest-neighbour divergence method based on
Rosenstein et al.~\cite{rosenstein1993practical}.
No additional delay embedding is applied.
For temporally separated neighbour pairs, let $d_j(i)$
denote their separation after $i$ sampling steps. We fit
\begin{equation}
  L(i)
  =
  \left\langle
    \log\frac{d_j(i)}{d_j(0)}
  \right\rangle_j
  \simeq b+\hat\lambda_1 i\Delta t,
  \label{eq:rosenstein}
\end{equation}
over a linear-growth interval, where $\Delta t$ is the
output sampling interval. A positive estimate supports
chaotic dynamics. We compare $\hat\lambda_1$ with the
reported NE6 value $\lambda_1\approx0.0280$
\cite{jafari2013elementary}.

\subsection{Conditional Lyapunov exponent}
The reported conditional Lyapunov exponent (CLE) measures
local response stability along the autonomous prediction.
We integrate the variational equation
\begin{equation}
  \dot{\xi}
  =
  D_{\bm{\theta}}R
  \bigl(\hat{\bm{\theta}}(t),\hat u(t)\bigr)\xi
\end{equation}
with periodic tangent-vector renormalization.
The generated signal $\hat u(t)$ is held fixed when taking
the derivative, so the feedback derivative is excluded.
A negative CLE indicates local contraction for this signal;
it does not establish global synchronization or the echo
state property for all inputs~\cite{lu2018attractor}.

\subsection{Normalized mean squared error (NMSE)}
Short-term accuracy is measured over the first
$M_{\mathrm p}=800$ prediction samples using the NMSE:
\begin{equation}
  \mathrm{NMSE}
  =
  \frac{1}{m}\sum_{\ell=1}^{m}
  \frac{
    \sum_{k=1}^{M_{\mathrm p}}
    (\hat u_{\ell,k}-u_{\ell,k})^2
  }{
    \sum_{k=1}^{M_{\mathrm p}}u_{\ell,k}^{2}
  }.
  \label{eq:nmse}
\end{equation}
This metric is evaluated in the scaled coordinates without
subtracting component means.
We classify $\mathrm{NMSE}\ge0.3$ as a failed short-term
prediction.

\subsection{Kullback--Leibler (KL) divergence}
Distributional agreement is quantified using common
three-dimensional histograms of the reference and predicted
trajectories. Both trajectories are standardized using the
reference mean and standard deviation, with $50$ bins per
coordinate. A pseudocount $\epsilon=10^{-12}$ is added to
each bin before normalization.

Writing $p_b$ and $q_b$ for the resulting reference and
prediction probabilities, respectively, we compute
\begin{equation}
  D_{\mathrm{KL}}(p\|q)
  =
  \sum_b p_b\log\frac{p_b}{q_b}.
  \label{eq:KL}
\end{equation}
Smaller values indicate closer empirical occupancy
distributions. Pseudocount smoothing keeps the estimate
finite even when the prediction misses reference-occupied
bins.

\subsection{Numerical assessment of hiddenness}
An attractor is hidden if its basin of attraction does not intersect
sufficiently small neighbourhoods of any equilibrium. For the all-to-all
and ring reservoirs, we first search for equilibria of the autonomous
system using multistart root finding, initialized from perturbed reservoir
states and equilibria found at neighbouring parameter values. A candidate
equilibrium is accepted when the infinity norm of the vector-field
residual is at most $10^{-8}$. At detected equilibria within the prescribed
output domain, we compute the Jacobian and integrate forward from small
perturbations along selected unstable eigendirections. We compare the
resulting long-time read-out distributions with that of the reconstructed
attractor using the KL divergence. If no equilibrium is found, or no
tested forward trajectory reproduces the attractor, we classify it as a
\emph{candidate hidden attractor}. This numerical classification does not
prove hiddenness, since the search may miss equilibria or trajectories
connecting their neighbourhoods to the attractor.

\section{Results}
\label{Results}

We compare the all-to-all, ring, and modified ring reservoirs over $420$
common parameter pairs. Throughout the figures, markers identify confined
chaotic reconstructions retained by the equilibrium and trajectory
screening. For the all-to-all and ring reservoirs, these are
\emph{candidate hidden attractors}: no equilibrium was found numerically,
or no tested continuation reproduced the attractor. These finite tests
do not prove hiddenness. For the modified ring, equilibria are excluded
by the governing equations. Crosses mark eligible reconstructions with
$\hat{\lambda}_1 \approx 0.03$ and
$D_{\mathrm{KL}}(p\,\|\,q)<0.25$; circles mark the remaining eligible
candidate hidden attractors.

\subsection{All-to-all topology}
\label{sec:alltoall}

\begin{figure}[tbp]
    \centering
    \includegraphics[width=\textwidth]{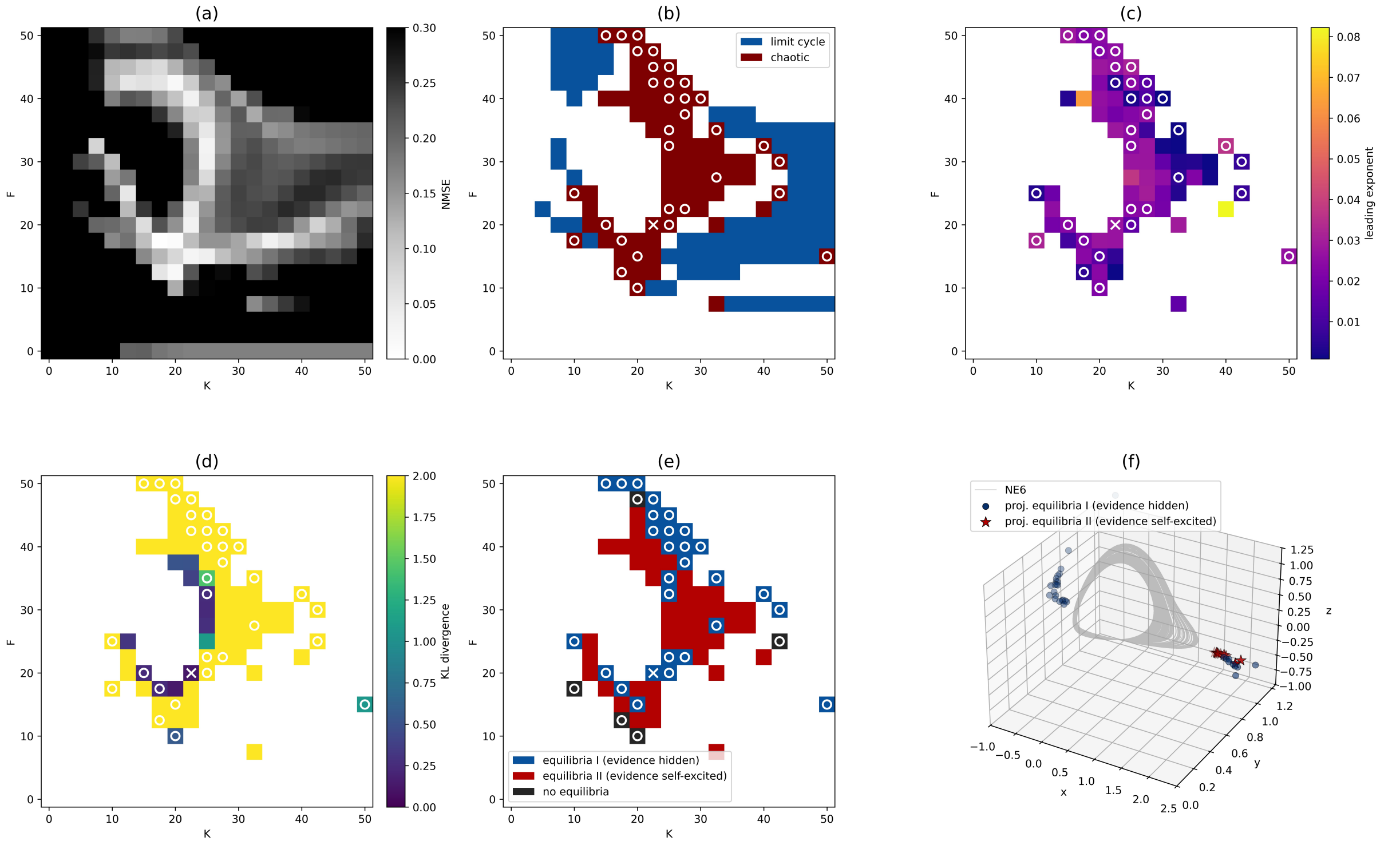}
    \caption{All-to-all reservoir trained on NE6.
    (a) NMSE. (b) Confined autonomous dynamics with negative maximal
    conditional Lyapunov exponent.
    (c) Leading Lyapunov exponent and (d) KL divergence for chaotic
    reconstructions.
    (e) Equilibrium-search and trajectory-screening outcomes.
    (f) Equilibria projected into output space alongside the attractor.
    Markers follow the convention defined at the beginning of this
    section.}
    \label{fig:naive_all_to_all}
\end{figure}

Accurate short-term prediction is confined to a wrench-shaped region of
the $(K,F)$-plane, with NMSE $\geq0.3$ at $223$ parameter pairs
(Fig.~\ref{fig:naive_all_to_all}a). Confined outputs with negative maximal
conditional Lyapunov exponent approach either limit cycles or chaotic
attractors; convergence to an equilibrium was not observed
(Fig.~\ref{fig:naive_all_to_all}b). Low NMSE does not ensure chaotic
reconstruction: the band $F\approx15$--$20$ at larger $K$ contains limit
cycles. Chaotic reconstructions have leading Lyapunov exponents mostly
between $0.01$ and $0.04$, comparable to $0.0280$ for NE6, but low KL
divergence occurs only in a narrow strip
(Figs.~\ref{fig:naive_all_to_all}c,d).

\begin{figure}[tbp]
    \centering
    \begin{subfigure}[t]{0.48\columnwidth}
        \centering
        \includegraphics[width=\linewidth]{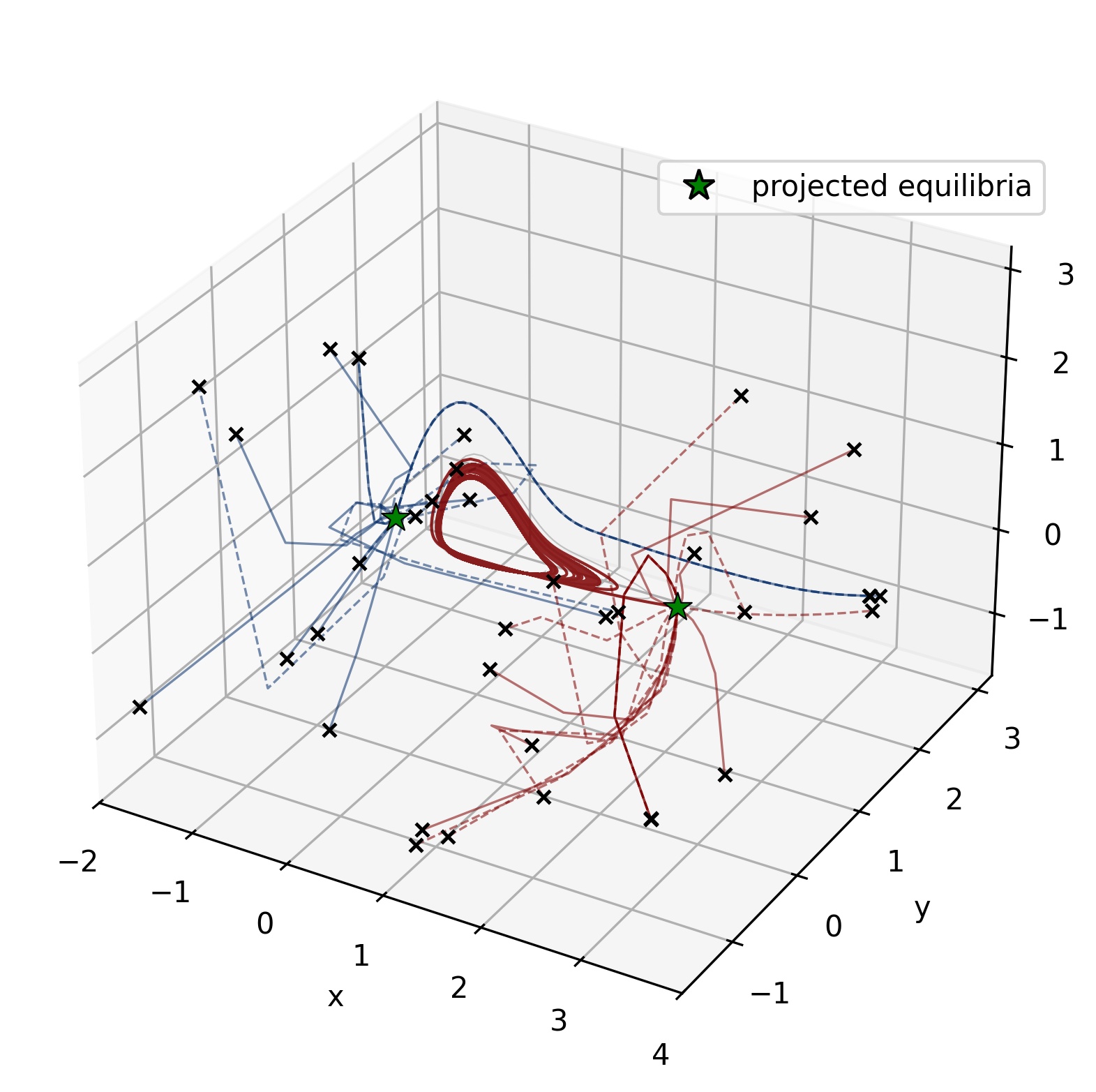}
        \caption{Matching continuation}
        \label{fig:SEA}
    \end{subfigure}
    \hfill
    \begin{subfigure}[t]{0.48\columnwidth}
        \centering
        \includegraphics[width=\linewidth]{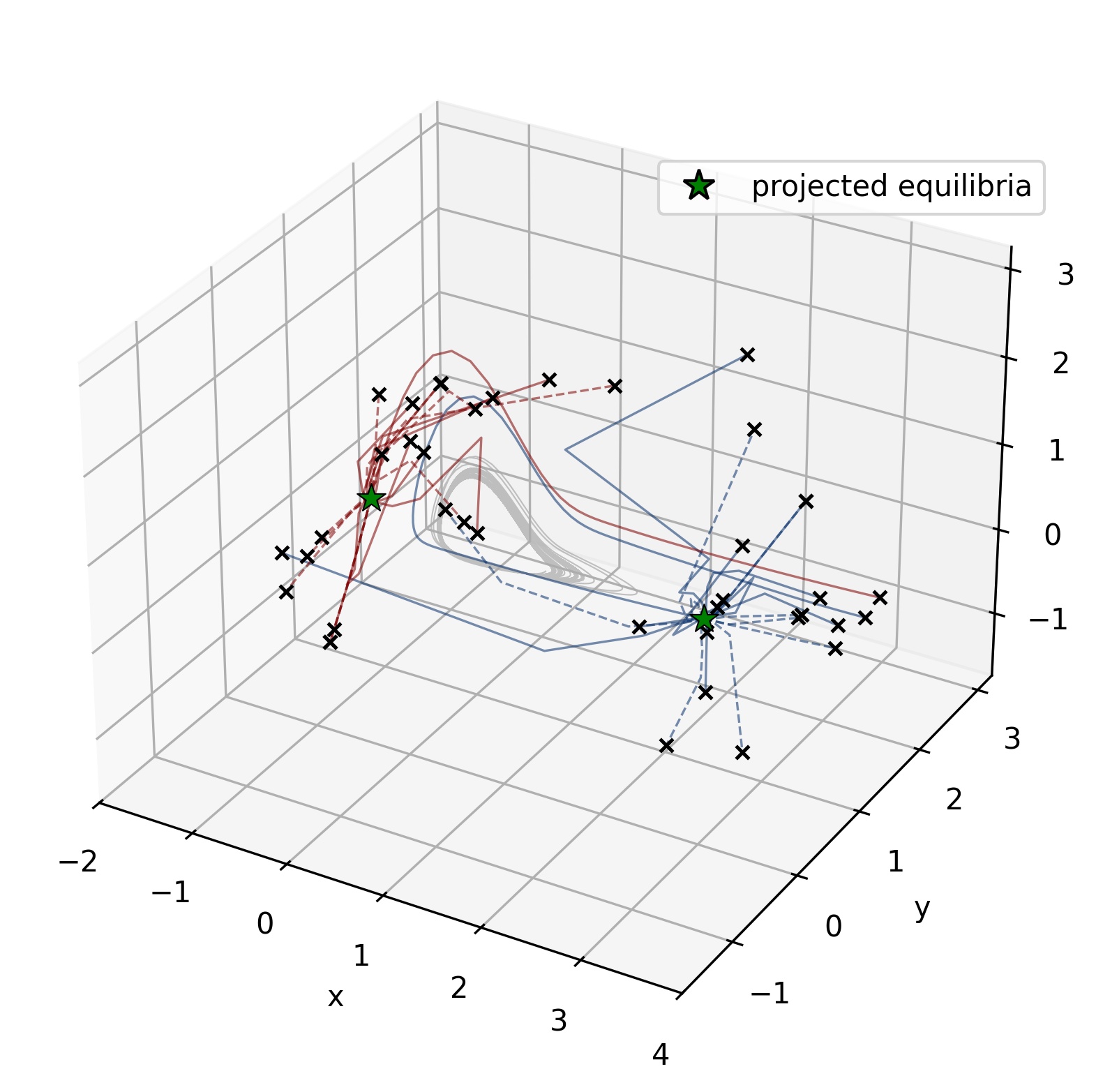}
        \caption{No matching continuation}
        \label{fig:EHA}
    \end{subfigure}
    \caption{Trajectories initialized near equilibria along selected
    Jacobian eigendirections at (a) $(K,F)=(30,22.5)$ and
    (b) $(K,F)=(30,40)$. Backward and forward trajectories test basin accessibility}
    \label{fig:continuation_examples}
\end{figure}

Unlike NE6, the learned reservoir can possess equilibria
(Figs.~\ref{fig:naive_all_to_all}e,f). Numerical screening retains $33$
candidate hidden attractors, only one of which meets both statistical
criteria, near $(K,F)=(22.5,20)$. Figure~\ref{fig:continuation_examples}
illustrates a matching continuation at $(30,22.5)$ and no matching
continuation at $(30,40)$. Only a matching forward trajectory supports
accessibility from an equilibrium neighbourhood; projected equilibrium
locations alone do not establish basin geometry.

\begin{figure}[tbp]
    \centering
    \includegraphics[width=\columnwidth]{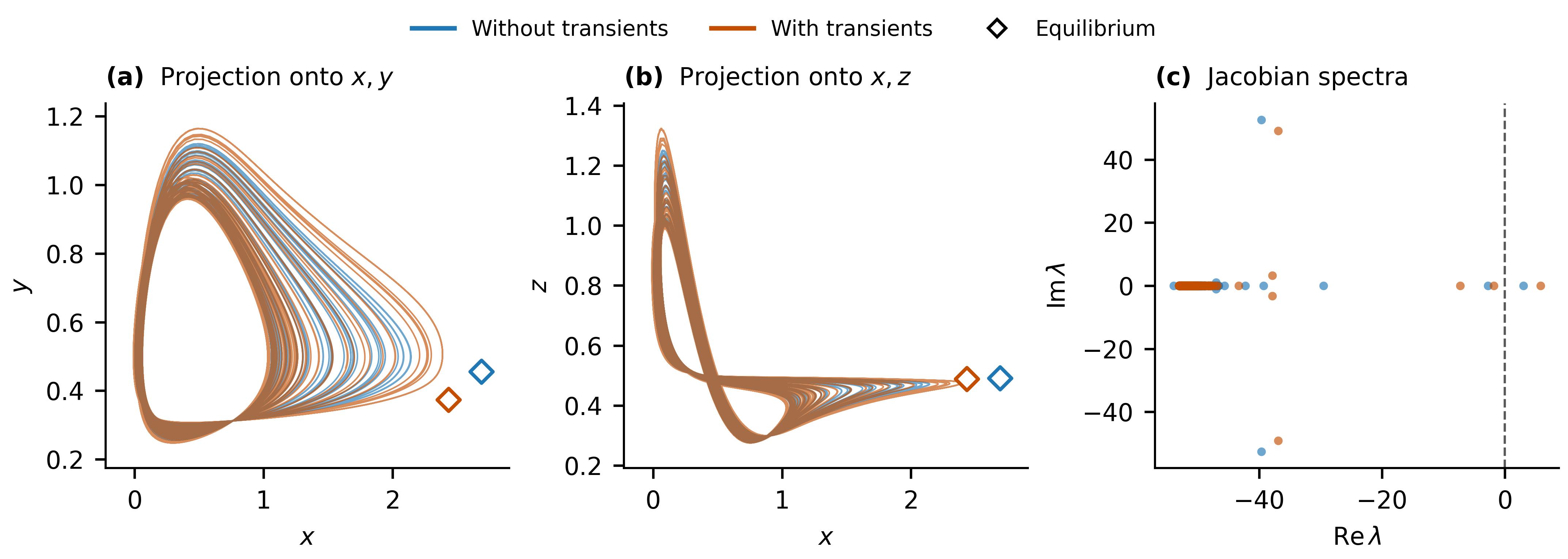}
    \caption{Training on attractor data alone (blue) and with additional
    transients (orange): reconstructed attractors projected onto
    (a) $(x,y)$ and (b) $(x,z)$.
    (c) Jacobian spectra at the detected equilibria, which persist and
    remain unstable in both cases.}
    \label{fig:perturbed_RC}
\end{figure}

Adding transient NE6 trajectories to training produces similar
reconstructed attractors but shifts rather than removes the detected
equilibrium (Fig.~\ref{fig:perturbed_RC}). Its Jacobian has a similar spectrum. Thus, additional training data do not remove the
equilibrium in this example, motivating changes to the reservoir
architecture.

\subsection{Ring topology}
\label{sec:RC_Ring}

\begin{figure}[tbp]
    \centering
    \includegraphics[width=\textwidth]{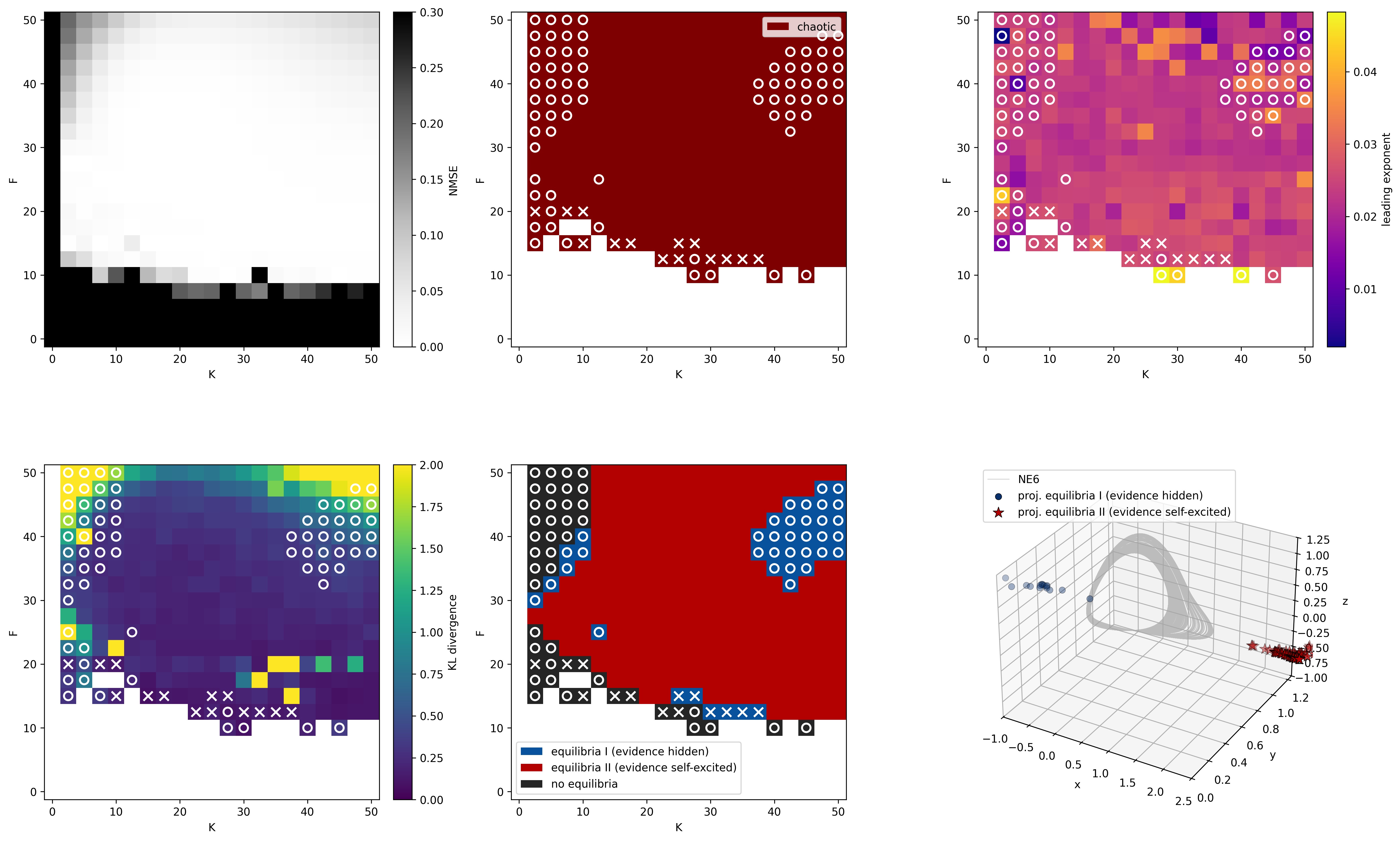}
    \caption{Ring reservoir trained on NE6, with the same panels and
    marker convention as Fig.~\ref{fig:naive_all_to_all}. The confined
    reconstructions with negative maximal conditional Lyapunov exponent
    shown in (b) are chaotic.}
    \label{fig:RC_Ring}
\end{figure}

Ring coupling substantially enlarges the low-NMSE region, approximately
over $F\gtrsim10$ and $K\gtrsim5$, reducing the number of parameter pairs
with NMSE $\geq0.3$ from $223$ to $75$
(Fig.~\ref{fig:RC_Ring}a). The confined region with negative maximal
conditional Lyapunov exponent is also larger, and the reconstructions
shown there are chaotic (Fig.~\ref{fig:RC_Ring}b). Leading Lyapunov
exponents mostly lie between $0.02$ and $0.03$, while low KL divergence
extends over a wider band than for all-to-all coupling
(Figs.~\ref{fig:RC_Ring}c,d).

Equilibria nevertheless persist throughout much of the accurate
reconstruction region (Figs.~\ref{fig:RC_Ring}e,f). Screening retains
$86$ candidate hidden attractors, including $14$ meeting both statistical
criteria. Ring coupling therefore improves reconstruction without
necessarily preserving NE6's absence of equilibria.

\subsection{Modified ring topology}
\label{sec:mod_Ring}
\begin{figure}[tbp]
    \centering
    \includegraphics[width=\textwidth]{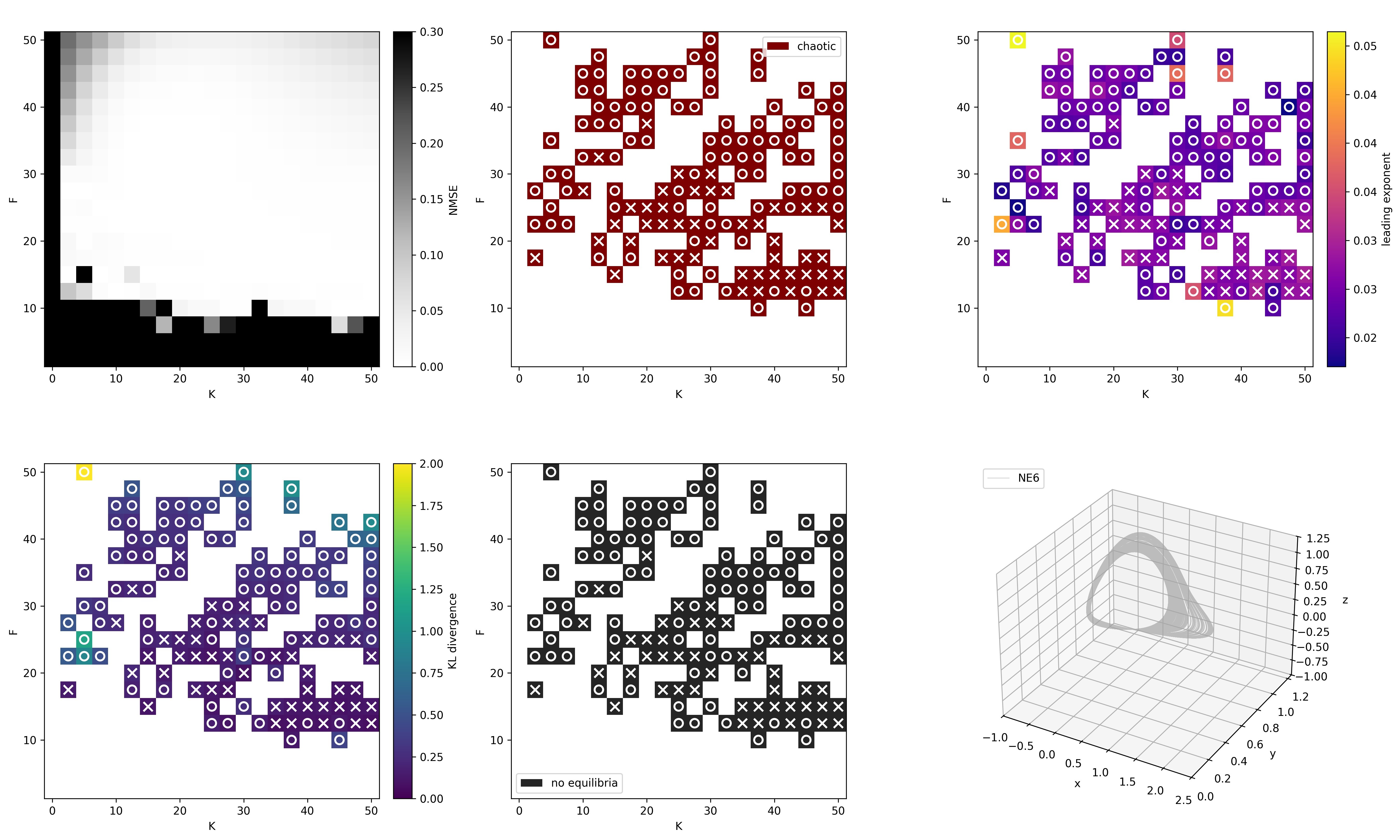}
    \caption{Modified ring reservoir trained on NE6, using the same
    layout and marker convention as Fig.~\ref{fig:naive_all_to_all}.
    Equilibria are excluded by the governing equations throughout the
    parameter range; marked attractors are therefore hidden.}
    \label{fig:RC_ModRing}
\end{figure}

The modified ring excludes equilibria by construction, so its attracting solutions are hidden under the equilibrium-based definition. We note that modified all-to-all yields poor performance as all nodes are connected to the modified equation. Hence, the results are not reported here. The ring network provides sparseness which ensures that the effect of the modified equation on the computational performance is limited. Short-term
accuracy remains comparable to the ordinary ring: $83$ parameter pairs
have NMSE $\geq0.3$ (Fig.~\ref{fig:RC_ModRing}a). Confined chaotic
reconstructions with negative maximal conditional Lyapunov exponent
occupy a broad region; their leading Lyapunov exponents mostly range
from $0.02$ to $0.035$, and low KL divergence occurs over a substantial
part of the parameter plane (Figs.~\ref{fig:RC_ModRing}b--d).
There are $144$ marked cases, of which $47$ meet both statistical criteria.

The modification adds only eight high-NMSE cases relative to the ring,
a difference of $1.9$ percentage points, while increasing the number of
crosses from $14$ to $47$ (Table~\ref{tab:system_comparison}). Among marked
cases, the fraction meeting both statistical criteria rises from $3\%$
for all-to-all coupling to $16\%$ for the ring and $33\%$ for the modified
ring. Thus, the ordinary ring has the fewest short-term prediction
failures, whereas the modified ring yields the most marked cases meeting
the statistical criteria and excludes equilibria analytically. These
counts describe the sampled parameter grid.

\begin{table}[tbp]
    \centering
    \caption{Comparison over $420$ common parameter pairs. Crosses and
    circles are disjoint classes. The NMSE column is counted independently
    and may overlap either marker class.}
    \label{tab:system_comparison}
    \begin{tabular}{lrrr}
        \toprule
        System & $\times$ & $\circ$ & NMSE $\geq 0.3$ \\
        \midrule
        All-to-all    &  1 & 32 & 223 \\
        Ring          & 14 & 72 &  75 \\
        Modified ring & 47 & 97 &  83 \\
        \bottomrule
    \end{tabular}
\end{table}

A reservoir without equilibria cannot have the Echo State Property (ESP) for
every input. To make this precise, let $\bm\theta(t)$ and
$\tilde{\bm\theta}(t)$ be two trajectories of the reservoir driven by the same
input $u$. We say that the reservoir has the \emph{global ESP for the input
$u$} if
\begin{equation}
  d\bigl(\bm\theta(t),\tilde{\bm\theta}(t)\bigr)\longrightarrow0
  \qquad(t\to\infty)
\end{equation}
for every pair of initial states,

where
$d(\bm\theta,\bm\eta)=\min_{\bm k\in\mathbb{Z}^N}
\|\bm\theta-\bm\eta+2\pi\bm k\|_2$ is the distance on the torus. It has the
global ESP for a class of admissible inputs
if this holds for every input in the class.

\begin{lemma}[Constant inputs and the global ESP]
\label{lem:modified_ring_no_esp}
For $F>0$ and $\omega_N\ne0$, the modified ring
\eqref{eq:modified_ring} fails the global echo state property
for every constant input $u(t)\equiv u_0$.
\end{lemma}

\begin{proof}
Fix $u_0$ and let $\Phi_t$ denote the resulting autonomous
flow on $\mathbb{T}^N$, with torus distance $d$.
Equation~\eqref{eq:modified_ring} excludes equilibria:
the first $N-1$ equations would require $b_1=b_{N-1}=0$,
leaving $\dot{\theta}_N=\omega_N\ne0$.

Suppose the global ESP holds, so every pair of initial
states is asymptotically forgotten. Applying this to
$\bm{\theta}_0$ and $\Phi_s(\bm{\theta}_0)$ gives, by autonomy,
\[
  d\bigl(\Phi_{t+s}(\bm{\theta}_0),
         \Phi_t(\bm{\theta}_0)\bigr)\longrightarrow0
  \qquad(t\to\infty)
\]
for every fixed $s\ge0$.
Compactness provides a sequence $t_n\to\infty$ such that
$\Phi_{t_n}(\bm{\theta}_0)\to\bm{\theta}_*$.
Continuity therefore yields
\[
  \Phi_s(\bm{\theta}_*)
  =\lim_{n\to\infty}\Phi_{t_n+s}(\bm{\theta}_0)
  =\bm{\theta}_*
  \qquad\text{for every }s\ge0.
\]
Thus $\bm{\theta}_*$ is an equilibrium, a contradiction.
\end{proof}

\section{Discussion}
\label{discussion}

We investigated whether driven Kuramoto reservoirs can reproduce both the
dynamics and the hidden character of the NE6 attractor. Ring coupling
broadened the parameter region supporting accurate short-term prediction
and agreement with the reference Lyapunov exponent and sampled
distribution. The modified ring retained comparable short-term accuracy
while yielding $47$ marked cases meeting the Lyapunov and KL criteria,
compared with $14$ for the ordinary ring and $1$ for all-to-all coupling
(Table~\ref{tab:system_comparison}).

Accurate output reconstruction does not determine hiddenness, which
depends on equilibria and their relationship to the basin of attraction
in the full reservoir state space. The ordinary reservoirs can acquire
equilibria absent from NE6, and finite equilibrium searches and trajectory
tests provide only candidate hidden attractors. Adding transient training
data did not remove the detected equilibrium in the example examined;
this does not establish a general limitation of transient training.
The modified ring instead excludes equilibria independently of the
trained readout (Sec.~\ref{sec:mod_Ring}). Consequently, any attracting
set realized by these dynamics is hidden, although its accuracy and
chaotic character require separate assessment. Excluding equilibria
preserves this structural property of NE6 without establishing agreement
of the remaining basin geometry.

The modified ring exhibits a more fragmented parameter region of
confined chaotic reconstructions than the ordinary ring
(Figs.~\ref{fig:RC_Ring}b and~\ref{fig:RC_ModRing}b).
If this fragmentation reflects genuine transitions between regular and
chaotic autonomous dynamics, the more intricate boundary could place
a larger fraction of sampled parameter pairs near the onset of chaos,
potentially contributing to successful reconstruction. This offers an
interpretation in terms of the edge-of-chaos hypothesis~\cite{carroll2020reservoir}.

The modification also limits claims about the echo state property.
For constant inputs, the absence of equilibria on the compact phase
torus prevents global forgetting of initial states. This is compatible
with local conditional contraction under a time-dependent input.
Accordingly, a negative maximal conditional Lyapunov exponent along the
NE6-driven response indicates local stability, rather than establishing
the global ESP. Stability claims should specify both the driving inputs
and the set of initial states considered.

These conclusions concern one target system, one reservoir family, and
the realizations and parameter ranges examined. Tests with additional
systems, reservoir realizations, and noisy data are needed to assess
generality. Systematic forward integrations from progressively smaller
equilibrium neighbourhoods would strengthen the basin analysis. A
further theoretical question is which nonconstant inputs and domains
of initial states support reliable synchronization in reservoirs whose
dynamics exclude equilibria.

\section*{Declaration of generative AI use}
ChatGPT was used for language editing of the manuscript and code support. The authors have thoroughly checked and revised the AI-assisted text and accept full responsibility for the article's content.
\section*{Declaration of competing interest}
The authors declare that they have no known competing financial interests or personal relationships that could have appeared to influence the work reported in this paper.

\section*{Acknowledgements}

The authors would like to thank Hebei Normal University for hosting Emel Khan. During this research, Thomas de Jong was also affiliated to Kanazawa University and Duke-Kunshan. This work was partially supported by JST CREST grant number JPMJCR2014.

\bibliographystyle{ieeetr}

\bibliography{cas-refs}

\end{document}